\documentclass[]{arxiv}

\title{Successive Convexification for 6-DoF Mars \\
       Rocket Powered Landing with Free-Final-Time}

\author{
  Michael Szmuk and Beh\c cet\ A\c c\i kme\c se \\
  {\normalsize\itshape University of Washington, Seattle, WA 98195-2400}
}

\usepackage{amsmath}
\usepackage{amsthm}
\usepackage{amssymb}
\usepackage{graphicx}
\usepackage{float}
\usepackage{hhline}
\usepackage{enumerate}
\usepackage{epstopdf}
\usepackage{multicol}
\usepackage{textcomp}       
\usepackage{gensymb}        
\usepackage{url}
\usepackage{nicefrac}

\newcommand{\bvec}[1]{\boldsymbol{\mathbf{#1}}}
\newcommand{\mat}[1]{\begin{bmatrix}#1\end{bmatrix}}
\newcommand{\definedas}{\triangleq}

\newcommand{\realD}[1]{\mathbb{R}^{#1}}

\newcommand{\realp}{\mathbb{R}_{+}}
\newcommand{\realpp}{\mathbb{R}_{++}}
\newcommand{\sthree}{\mathcal{S}^3}
\newcommand{\sothree}{SO(3)}
\newcommand{\sympd}[1]{\mathbb{S}^{#1}_{++}}
\newcommand{\onenorm}[1]{\|#1\|_1}
\newcommand{\twonorm}[1]{\|#1\|_2}

\newcommand{\inertial}{\mathcal{I}}
\newcommand{\body}{\mathcal{B}}
\newcommand{\refframe}[1]{$\mathcal{F}_{\mathcal{#1}}$}

\newcommand{\skewmat}[1]{\left[#1\times\right]}
\newcommand{\OMEGA}[1]{\Omega\big(#1\big)}
\newcommand{\gI}{\bvec{g}_\inertial}
\newcommand{\JB}{J_\body}
\newcommand{\rTB}{\bvec{r}_{T,\body}}
\newcommand{\alphamdot}{\alpha_{\dot{m}}}
\newcommand{\omegamax}{\omega_{max}}
\newcommand{\maxgimbalangle}{\delta_{max}}
\newcommand{\tilt}{\theta}
\newcommand{\maxtilt}{\theta_{max}}
\newcommand{\Tmin}{T_{min}}
\newcommand{\Tmax}{T_{max}}
\newcommand{\glideslope}{\gamma_{gs}}
\newcommand{\cIB}{C_{\body/\inertial}}
\newcommand{\cBI}{C_{\inertial/\body}}
\newcommand{\Isp}{I_{sp}}
\newcommand{\gs}{g_0}
\newcommand{\FI}{\bvec{F}_\inertial}
\newcommand{\MB}{\bvec{M}_\body}
\newcommand{\eu}{\bvec{e}_1}
\newcommand{\ee}{\bvec{e}_2}
\newcommand{\en}{\bvec{e}_3}

\newcommand{\m}{m}
\newcommand{\rI}{\bvec{r}_\inertial}
\newcommand{\vI}{\bvec{v}_\inertial}
\newcommand{\qIB}{q_{\body/\inertial}}
\newcommand{\omegaB}{\bvec{\omega}_\body}

\newcommand{\TB}{\bvec{T}_\body}

\newcommand{\mdot}{\dot{m}}
\newcommand{\rIdot}{\dot{\bvec{r}}_\inertial}
\newcommand{\vIdot}{\dot{\bvec{v}}_\inertial}
\newcommand{\qIBdot}{\dot{q}_{\body/\inertial}}
\newcommand{\omegaBdot}{\dot{\bvec{\omega}}_\body}

\newcommand{\mi}{\m_{wet}}
\newcommand{\rIi}{\bvec{r}_{\inertial,i}}
\newcommand{\vIi}{\bvec{v}_{\inertial,i}}

\newcommand{\omegaBi}{\bvec{\omega}_{\body,i}}

\newcommand{\mf}{\m_{dry}}
\newcommand{\rIf}{\bvec{0}}
\newcommand{\vIf}{\bvec{v}_{\inertial,f}}
\newcommand{\qIBf}{q_{\body/\inertial,f}}
\newcommand{\omegaBf}{\bvec{0}}

\newcommand{\mki}{\m_{0}}
\newcommand{\rIki}{\bvec{r}_{\inertial,0}}
\newcommand{\vIki}{\bvec{v}_{\inertial,0}}

\newcommand{\omegaBki}{\bvec{\omega}_{\body,0}}

\newcommand{\mk}{\m_{k}}
\newcommand{\rIk}{\bvec{r}_{\inertial,k}}
\newcommand{\vIk}{\bvec{v}_{\inertial,k}}
\newcommand{\qIBk}{q_{\body/\inertial,k}}
\newcommand{\omegaBk}{\bvec{\omega}_{\body,k}}
\newcommand{\nuk}{\bvec{\nu}_k}

\newcommand{\rIkf}{\bvec{r}_{\inertial,K}}
\newcommand{\vIkf}{\bvec{v}_{\inertial,K}}
\newcommand{\qIBkf}{q_{\body/\inertial,K}}
\newcommand{\omegaBkf}{\bvec{\omega}_{\body,K}}
\newcommand{\uukf}{\bvec{u}_{K}}

\newcommand{\ti}{0}
\newcommand{\tf}{t_f}
\newcommand{\tauk}{\tau_k}
\newcommand{\taukp}{\tau_{k+1}}

\newcommand{\dd}[2]{\frac{d#1}{d#2}}
\newcommand{\xx}{\bvec{x}}
\newcommand{\uu}{\bvec{u}}
\newcommand{\zz}{\bvec{z}}
\newcommand{\xxo}{\hat{\xx}}
\newcommand{\uuo}{\hat{\uu}}

\newcommand{\sso}{\hat{\sigma}}

\newcommand{\xxk}{\xx_k}
\newcommand{\uuk}{\uu_k}
\newcommand{\xxkp}{\xx_{k+1}}
\newcommand{\uukp}{\uu_{k+1}}
\newcommand{\lambdal}{\alpha_k}
\newcommand{\lambdar}{\beta_k}
\newcommand{\stm}[2]{\Phi_{A}\big(#1,#2\big)}
\newcommand{\AD}{\bar{A}_k}
\newcommand{\BD}{\bar{B}_k}
\newcommand{\CD}{\bar{C}_k}
\newcommand{\SD}{\bar{\Sigma}_k}
\newcommand{\zD}{\bar{\bvec{z}}_k}

\newcommand{\f}[2]{f\big(#1,#2\big)}
\newcommand{\dfdx}{\frac{\partial}{\partial\xx}f(\xx,\uu)\Bigr|_{\xxo(\tau),\uuo(\tau)}}
\newcommand{\dfdu}{\frac{\partial}{\partial\uu}f(\xx,\uu)\Bigr|_{\xxo(\tau),\uuo(\tau)}}

\newcommand{\g}[1]{g\big(#1\big)}

\newcommand{\wnu}{w_{\nu}}
\newcommand{\wdelta}{w_{\Delta}^i}
\newcommand{\wdeltasigma}{w_{\Delta_{\sigma}}}

\newcommand{\nubar}{\bar{\bvec{\nu}}}
\newcommand{\Deltabar}{\bar{\bvec{\Delta}}}

\newcommand{\setK}{\mathcal{K}}
\newcommand{\setKm}{\bar{\mathcal{K}}}

\newcommand{\UL}{U_L}
\newcommand{\UM}{U_M}
\newcommand{\UT}{U_T}

\newcommand{\boxing}[6]{
	\begin{figure*}[#1]
		\noindent\makebox[\textwidth][c] {
			\fbox{
				\begin{minipage}{#4}
					\vspace{0.2cm}
					\begin{#2} {\bf : \emph{#5}} \label{#3} \end{#2}
					#6
				\end{minipage}
			}
		}
	\end{figure*}
}

\begin{document}
\maketitle


\begin{abstract}

In this paper, we employ successive convexification to solve the minimum-time 6-DoF rocket powered landing problem. The contribution of this paper is the development and demonstration of a free-final-time problem formulation that can be solved iteratively using a successive convexification framework. This paper is an extension of our previous work on the 3-DoF \emph{free-final-time} and the 6-DoF \emph{fixed-final-time} minimum-fuel problems. Herein, the vehicle is modeled as a 6-DoF rigid-body controlled by a single gimbaled rocket engine. The trajectory is subject to a variety of convex and non-convex state- and control-constraints, and aerodynamic effects are assumed negligible. The objective of the problem is to determine the optimal \emph{thrust commands} that will minimize the time-of-flight while satisfying the aforementioned constraints. Solving this problem quickly and reliably is challenging because (a) it is nonlinear and non-convex, (b) the validity of the solution is heavily dependent on the accuracy of the discretization scheme, and (c) it can be difficult to select a suitable reference trajectory to initialize an iterative solution process. To deal with these issues, our algorithm (a) uses successive convexification to eliminate non-convexities, (b) computes the discrete linear-time-variant system matrices to ensure that the converged solution perfectly satisfies the original nonlinear dynamics, and (c) can be initialized with a simple, dynamically inconsistent reference trajectory. Using the proposed convex formulation and successive convexification framework, we are able to convert the original non-convex problem into a sequence of convex second-order cone programming (SOCP) sub-problems. Through the use of Interior Point Method (IPM) solvers, this sequence can be solved quickly and reliably, thus enabling higher fidelity real-time guidance for rocket powered landings on Mars.

\end{abstract}


\section*{Nomenclature}
\begin{multicols}{2}
	\begin{tabbing}
	\hspace*{1.6cm}\= \kill
	DoF \> Degrees-of-Freedom \\ 
	IPM \> Interior Point Method \\
	SC \> Successive Convexification \\
	SOCP \> Second-Order Cone Programming \\
	DCM \> Direction cosine matrix \\
	$\bvec{e}_i$ \> Unit vector pointing along $i^{th}$-axis \\
	$\inertial$ \> Subscript used to denote the inertial frame \\
	$\body$ \> Subscript used to denote the body frame \\
	\refframe{I} \> The inertial NED frame \\
	\refframe{B} \> The body frame \\
	$t$ \> Time \\
	$\tau$ \> Normalized trajectory time \\
	$\tf$ \> Time-of-flight \\
	$\mi$ \> Dry mass of the landing vehicle \\
	$\mf$ \> Wet mass of the landing vehicle \\
	$\rTB$ \> Gimbal-point position vector \\ 
	$\gI$ \> Gravity vector expressed in \refframe{I} \\
	$m$ \> Vehicle mass \\
	$\rI$ \> Inertial position of the vehicle \\
	$\vI$ \> Inertial velocity of the vehicle \\
	$\qIB$ \> Unit quaternion \\
	$\omegaB$ \> Angular velocity vector \\
	$\FI$ \> Force acting on the vehicle \\
	$\TB$ \> Commanded thrust vector \\
	$\bvec{\nu}$ \> Virtual control vector \\
	$\cIB$ \> DCM corresponding to $\qIB{}$ \\
	$\tilt$ \> Vehicle tilt angle \\
	$\delta$ \> Gimbal angle \\
	$\JB$ \> Moment of inertia of the vehicle \\
	$\glideslope$ \> Glide-slope cone constraint angle \\
	$\maxtilt$ \> Maximum allowable tilt angle \\
	$\omegamax$ \> Maximum allowable angular rate \\
	$\maxgimbalangle$ \> Maximum allowable gimbal angle \\
	\end{tabbing}
\end{multicols}


\clearpage
\section{Introduction}

In this paper, we employ successive convexification to solve the minimum-time 6-DoF rocket powered landing problem. The work presented herein builds on our previous work on the 3-DoF free-final-time\cite{szmuk2016successive} and the 6-DoF fixed-final-time\cite{szmuk2017successive} minimum-fuel problems. The contribution of this paper is the development and demonstration of a 6-DoF free-final-time formulation that can be used within a successive convexification framework to solve the non-convex minimum-time 6-DoF rocket landing problem.

Autonomous optimal guidance technology that enables precision rocket powered landing is important for one primary reason: it provides a way to recover from a wider range of dispersions and uncertainties encountered during the entry, descent, and landing phase. This increases the probability that the vehicle will land safely and at a desirable location, and thus increases the likelihood of mission success. Such capabilities may be used to improve the scientific return of planetary science missions through the selection of more challenging and scientifically interesting landing sites, or by requiring less propellant mass to achieve the same level of landing robustness. This technology also enables the robust recovery of vertical-takeoff-vertical-landing (VTVL) launch vehicles\cite{larsNAE}, thus allowing for dramatic launch cost reductions through the reuse of boosters. The importance of the latter point is highlighted by the efforts and triumphs of commercial space companies to land and reuse launch vehicles in recent years.

Solving autonomous optimal landing guidance problems quickly and reliably is challenging for three reasons. First, they are nonlinear and non-convex. Second, for numerical implementations, the validity of the solution is heavily dependent on the accuracy of the discretization scheme. Third, it can be difficult to select a suitable reference trajectory to initialize an iterative solution process.

Several insightful methods have been proposed \cite{klumpp,najson_pdg,topcu_pdg} to analytically obtain sub-optimal guidance trajectory solutions. Although some methods successfully incorporated control constraints, attitude, and attitude-rate into the problem\cite{sostaric2005powered}, these analytical solutions could only handle a limited set of typical mission constraints.

In contrast, convex optimization has been considered a prime candidate for on-board autonomous guidance applications due to its deterministic behavior, its global optimality guarantees\cite{BoydConvex}, its ability to provide certificates of convergence and infeasibility, and the availability of efficient interior point method (IPM) algorithms\cite{dueri_ifac14_ipm,dueri2016customized,cvx,domahidi2012efficient} to solve convex problems. However, as mentioned before, the primary roadblock to employing convex optimization in this context is the non-convex nature of the problem.

Lossless convexification was one method proposed to circumvent the non-convex minimum thrust constraint of the landing problem.\cite{behcet_aut11,matt_aut1,larssys12,behcetaaiaa05,behcetjgcd07,larsjgcd10,pointing2013,harris2013maximum} This method was successfully demonstrated numerous times during flight experiments\cite{gfold,gfold_ieee_aero_14}. In a different instance, a dual-quaternion-based approach was proposed in which a 6-DoF line-of-sight constraint was convexified.\cite{mesbahi_6dof,lee2016constrained} Since this method was based on dual-quaternions, it was inherently equipped to handle 6-DoF motion, but relied on piecewise-affine approximations of the nonlinear dynamics. As such, the solution degraded in accuracy as the temporal resolution of the discretization was reduced. Both of these methods were limited to fixed-final-time problems.

Although these convexification techniques proved to greatly simplify the solution process for certain non-convex problems, they remain applicable only to a narrow class of problems. On the other hand, successive convexification techniques offer a framework for solving more general non-convex optimal control problems, at the expense of more computational complexity and weakened optimality and convergence guarantees. These methods work by transforming the non-convex problem into a sequence of convex optimization problems that are solved in succession to convergence.\cite{jordithesis,liu2014solving,liu2015entry,SCvx_2016arXiv,SCvx_2017arXiv,SCvx_cdc16,wang2016constrained,liu_new}

In this paper, our strategy is to employ successive convexification to solve the minimum-time 6-DoF landing problem. First, we initialize the process using a simple, dynamically inconsistent reference trajectory. Then, we linearize and discretize the problem to generate a convex Second-Order Cone Programming (SOCP) sub-problem. We solve these sub-problems in succession, with each iteration linearizing around the previously generated solution. Lastly, we structure this process in a way that ensures the converged solution will satisfy the original nonlinear dynamics, kinematics, and constraints. We emphasize that our method does \emph{not} perform a line search on the time-of-flight. Furthermore, in contrast to existing heuristic strategies, our method is capable of generating trajectories that are dynamically feasible, that adhere to the prescribed constraints, and that are more optimal, thus enlarging the usable flight envelope. Thus, we believe that the proposed method is a promising option for future autonomous applications.

The paper is organized as follows. In Section~\ref{sec:prob_des}, we state the continuous-time non-convex problem description. In Section~\ref{sec:convex}, we describe the successive convexification framework, and propose a method for converting the original problem into a sequence of discrete-time convex optimization problems. In Section~\ref{sec:results}, we present simulation of the proposed algorithm. Lastly, in Section~\ref{sec:conclusion} we provide concluding remarks.


\clearpage
\section{Problem Description} \label{sec:prob_des}

In this section, we outline the 6-DoF minimum-time rocket powered landing problem in its original continuous-time non-convex form. This problem description is similar to its fixed-final-time minimum-fuel counterpart provided in our previous work\cite{szmuk2017successive}, with the notable exception that the final time of the trajectory is now an optimization variable. For completeness, we include the entirety of the problem description. In what ensues, we use $\tf$ to denote the final time of the trajectory, and refer to $\tf$ as the time-of-flight. Furthermore, we use the following two reference frames:
\begin{enumerate}
	\item \refframe{I}: An inertially-fixed Up-East-North reference frame with its origin located at the landing site.
	\item \refframe{B}: A body-fixed frame centered at the vehicle center-of-mass, with its X-axis pointing along the vertical axis of the vehicle (i.e. parallel to the thrust vector when the engine gimbal angle is zero), its Y-axis pointing out the side of the vehicle, and its Z-axis completing the right-handed system.
\end{enumerate}
\noindent The remainder of this section covers the dynamics and kinematics that govern the problem, the state and control constraints imposed on the trajectory, and the boundary conditions enforced at the beginning and end of the trajectory. The section concludes with a concise summary of the problem description.

\subsection{Dynamics and Kinematics}

Herein, we treat the vehicle as a rigid body subject to constant gravitational acceleration, $\gI\in\realD{3}$, and negligible aerodynamic forces. The vehicle is assumed to actuate a single gimbaled rocket engine to generate a thrust vector within a feasible range of magnitudes and gimbal angles. Note that with minor modifications, this formulation can be adapted for applications that use multiple non-gimbaled rocket motors. Furthermore, we assume the vehicle depletes its mass, $\m(t)\in\realpp$, at a rate proportional to the magnitude of the commanded thrust vector\cite{behcetjgcd07}, $\TB(t)\in\realD{3}$, which is expressed in \refframe{B} coordinates. For tractability, we assume that the inertia tensor and the position of the center-of-mass are constant despite the depletion of mass. The proportionality constant, $\alphamdot$, is given in terms of the vacuum-specific-impulse, $\Isp$, and Earth's standard gravity constant, $\gs$, as follows:
\begin{equation*}
	\alphamdot \definedas \frac{1}{\Isp\gs}
\end{equation*}
Thus, the mass depletion dynamics are given by:
\begin{equation} \label{prob_des:mass_depletion}
	\mdot(t) = -\alphamdot\twonorm{TB(t)}
\end{equation}
We express the position, velocity, and force acting on the vehicle in \refframe{I} coordinates, and denote them by $\rI(t)\in\realD{3}$, $\vI(t)\in\realD{3}$, and $\FI(t)\in\realD{3}$, respectively. Thus, the translational dynamics are given by:
\begin{subequations} \label{prob_des:trans_dyn}
\begin{align}
	\rIdot(t) &= \vI(t) \\
	\vIdot(t) &= \frac{1}{m(t)}\FI(t)+\gI
\end{align}
\end{subequations}
\noindent We use unit quaternions to parametrize the attitude of \refframe{B} relative to \refframe{I}, and denote them by $\qIB(t)\in\sthree$. In this paper, we use the leading scalar element convention for quaternions:
\begin{equation*}
	\qIB \definedas \mat{\cos(\xi/2) \\ \sin(\xi/2)\hat{\bvec{n}}} = \mat{q_0 & q_1 & q_2 & q_3}^T
\end{equation*}
\noindent The direction cosine matrix (DCM) that encodes the attitude transformation from \refframe{I} to \refframe{B} is denoted by $\cIB(t)\in\sothree$, where $\cIB$ is related to $\qIB$ through the following relation:
\begin{equation*}
	\cIB = \mat{ 1-2 (  q_2^2 + q_3^2  ) &   2 (q_1 q_2 + q_0 q_3) &   2 (q_1 q_3 - q_0 q_2) \\
	               2 (q_1 q_2 - q_0 q_3) & 1-2 (  q_1^2 + q_3^2  ) &   2 (q_2 q_3 + q_0 q_1) \\
			 	 			   2 (q_1 q_3 + q_0 q_2) &   2 (q_2 q_3 - q_0 q_1) & 1-2 (  q_1^2 + q_2^2  )}
\end{equation*}
The inverse transformation (i.e. from \refframe{B} to \refframe{I}) is denoted by $\cBI(t) \definedas \cIB^{-1}(t) = \cIB^T(t)$. We use $\omegaB(t)\in\realD{3}$ to denote the angular velocity vector of \refframe{B} relative to \refframe{I}, expressed in \refframe{B} coordinates, and, for some $\bvec{\xi}\in\realD{3}$, define the skew-symmetric matrices $\skewmat{\bvec{\xi}}$ and $\OMEGA{\bvec{\xi}}$ as:
\begin{equation*}
	\skewmat{\bvec{\xi}} \definedas \mat{0&-\xi_z&\xi_y\\ \xi_z&0&-\xi_x\\ -\xi_y&\xi_x&0} \quad\quad \OMEGA{\bvec{\xi}} \definedas \mat{0&-\xi_x&-\xi_y&-\xi_z\\ \xi_x &0&\xi_z&-\xi_y\\ \xi_y &-\xi_z&0&\xi_x\\ \xi_z&\xi_y&-\xi_x&0}
\end{equation*}
\noindent We denote the torque acting on the vehicle as $\MB(t)\in\realD{3}$, and the inertia tensor of the vehicle in \refframe{B} coordinates as $\JB\in\sympd{3}$. Thus, the quaternion kinematics and attitude dynamics are given by:
\begin{subequations} \label{prob_des:att_kin_dyn}
\begin{align} 
	\qIBdot(t) &= \frac{1}{2}\OMEGA{\omegaB(t)}\qIB(t) \\
	\JB\omegaBdot(t) &= \MB(t) - \skewmat{\omegaB(t)}\JB\omegaB(t)
\end{align}
\end{subequations}
\noindent Lastly, we must relate the force and torque vectors, $\FI(t)$ and $\MB(t)$, to the thrust vector, $\TB(t)$. The force is related to the thrust vector through the following relation:
\begin{equation} \label{prob_des:force}
	\FI(t) = \cBI(t)\TB(t)
\end{equation}
\noindent Based on our earlier assumption that the center-of-mass does not move, it follows that the moment arm from the center-of-mass to the gimbal point of the engine is constant. We denote this constant position vector in \refframe{B} coordinates as $\rTB\in\realD{3}$. Thus, the torque is given by:
\begin{equation} \label{prob_des:torque}
	\MB(t) = \skewmat{\rTB}\TB(t)
\end{equation}

\subsection{State and Control Constraints} \label{sec:state_control_constraints}

To begin, we restrict the mass of the vehicle to remain above the dry mass, $\mf$, using the following convex constraint:
\begin{equation} \label{prob_des:mdry}
	\mf \leq m(t)
\end{equation}
\noindent The path of the vehicle is restricted to lie within an upward-facing glide-slope cone that makes an angle $\glideslope\in[0\degree,90\degree)$ with the horizontal, and that is centered at the origin of \refframe{I}. This amounts to the convex constraint given by:
\begin{subequations} \label{prob_des:gs_cone}
	\begin{align}
		\eu\cdot\rI(t) & \geq \tan{\glideslope}\twonorm{H_{23}^T \rI(t)} \\
		H_{23} &\definedas \mat{\ee & \en}
	\end{align}
\end{subequations}
\noindent We define the tilt angle of the vehicle as the angle between the X-axes of \refframe{B} and \refframe{I}, and denote it by $\tilt(t)$. More formally, this is expressed as:
\begin{equation} \label{prob_des:tilt_expr}
	\cos{\tilt(t)} = \eu\cdot\cBI(t)\eu = 1-2\big(q_2^2(t)+q_3^2(t)\big)
\end{equation}
\noindent Note that the quaternion elements in Eq.~\ref{prob_des:tilt_expr} are elements of $\qIB(t)$. To avoid excessive tilt angles in the trajectory, we limit $\tilt(t)$ to a maximum value of $\maxtilt$. If we immerse $\qIB$ in $\realD{4}$, we can impose this limit through the following convex constraint:
\begin{equation} \label{prob_des:tilt}
	\cos{\maxtilt} \leq 1-2\big(q_2^2(t)+q_3^2(t)\big)
\end{equation}
\noindent Furthermore, we limit the vehicle to a maximum angular rate of $\omegamax$ using the following convex constraint:
\begin{equation} \label{prob_des:ang_rate}
	\twonorm{\omegaB(t)} \leq \omegamax
\end{equation}
\noindent Lastly, the commanded thrust vector is constrained to lie within the magnitude interval $[\Tmin,\Tmax]$ and the gimbal angle interval $[0,\maxgimbalangle]$, where $\maxgimbalangle\in(0\degree,90\degree)$:
\begin{align} 
	0 < \Tmin \leq \twonorm{\TB(t)} &\leq \Tmax \label{prob_des:thrust_mag} \\
	\cos{\maxgimbalangle}\twonorm{\TB(t)} &\leq \eu\cdot\TB(t) \label{prob_des:gimbal_ang}
\end{align}
\noindent Note that the upper thrust magnitude bound is a convex constraint, whereas the lower thrust magnitude bound is non-convex. For the range of $\maxgimbalangle$ considered, the gimbal angle constraint is convex.

\subsection{Boundary Conditions} \label{sec:prob_des_bcs}

Depending on the application, various combinations of boundary conditions may be enforced. In the scenario considered here, we enforce one common combination of boundary conditions, with the understanding that a different valid combination could be enforced with ease. First, we constrain the initial mass, position, velocity, and angular rate to $\mi$, $\rIi$, $\vIi$, and $\omegaBi$, respectively, thus leaving the initial attitude free. We constrain the final position, velocity, attitude, and angular rate to $\rIf$, $\vIf$, $\qIBf$, and $\omegaBf$, respectively, thus leaving the final mass free. Lastly, we leave the initial thrust vector unconstrained, while constraining the final thrust vector to point along the X-axis of \refframe{B} in order to null out the torques at the final condition.

\subsection{Problem Statement}
	
We conclude by providing the complete problem description for the minimum-time 6-DoF rocket powered landing problem. We emphasize that in its original form, the problem is a non-convex continuous-time free-final-time problem. The objective is to minimize the time-of-flight subject to (a) the boundary conditions described in Section \ref{sec:prob_des_bcs}, (b) the dynamics and kinematics in Eqs.~\ref{prob_des:mass_depletion}-\ref{prob_des:torque}, (c) the state constraints in Eqs.~\ref{prob_des:mdry}-\ref{prob_des:gs_cone} and \ref{prob_des:tilt}-\ref{prob_des:ang_rate}, and (d) the control constraints in Eqs.~\ref{prob_des:thrust_mag}-\ref{prob_des:gimbal_ang}. The free variables are the commanded thrust vector, $\TB(t)$, and the time-of-flight, $\tf$. Problem~\ref{problem:non_convex} summarizes the complete non-convex problem description.
\vspace{0.2cm}
\boxing{ht!}{problem}{problem:non_convex}{16cm}{Non-Convex Continuous-Time Free-Final-Time Problem}{
	\begin{tabular}{lll}
		&& \\
		\underline{Cost Function}: &\null\quad\quad\null& \\
		&& $\underset{\tf,\TB(t)}{\text{minimize}}\quad \tf$ \\
		&& \\
		&& subject to: \\
		\underline{Boundary Conditions}: && \\
		&& $ {\begin{aligned}
					\m(\ti) &= \mi && \\ 
					\rI(\ti) &= \rIi & \rI(\tf) &= \rIf \\
					\vI(\ti) &= \vIi & \vI(\tf) &= \vIf \\
					&& \qIB(\tf) &= \qIBf \\
					\omegaB(\ti) &= \omegaBi & \omegaB(\tf) &= \omegaBf \\
					&& \ee\cdot\TB(\tf) &= \en\cdot\TB(\tf) = 0
				\end{aligned}} $ \\
		\underline{Dynamics}: && \\
		&& $ {\begin{aligned}
					\dot{m}(t) &= -\alphamdot\twonorm{\TB(t)} \\
					\rIdot(t) &= \vI(t) \\
					\vIdot(t) &= \frac{1}{m(t)}\cBI(t)\TB(t)+\gI \\
					\qIBdot(t) &= \frac{1}{2}\OMEGA{\omegaB(t)}\qIB(t) \\
					\JB\omegaBdot(t) &= \skewmat{\rTB}\TB(t) - \skewmat{\omegaB(t)}\JB\omegaB(t)
				\end{aligned}} $ \\
		\underline{State Constraints}: && \\
		&& $ {\begin{aligned}
					\mf &\leq m(t) \\
					\tan{\glideslope}\twonorm{H_{23}\rI(t)} &\leq \eu\cdot\rI(t) \\
					\cos{\maxtilt} &\leq 1-2\big(q_2^2(t)+q_3^2(t)\big) \\
					\twonorm{\omegaB(t)} &\leq \omegamax \\
				\end{aligned}} $ \\
		\underline{Control Constraints}: && \\
		&& $ {\begin{aligned}
					0 < \Tmin \leq \twonorm{\TB(t)} &\leq \Tmax \\
					\cos{\maxgimbalangle}\twonorm{\TB(t)} &\leq \eu\cdot\TB(t)
				\end{aligned}} $ \\
		&& \\
	\end{tabular}
}


\clearpage
\section{Convex Formulation} \label{sec:convex}

In this section, we develop the convex formulation for Problem~\ref{problem:non_convex}. Through a series of modifications, we morph the non-convex free-final-time problem into a convex fixed-final-time sub-problem, specifically an SOCP problem. This convex sub-problem will ultimately be solved repeatedly to convergence, a process we refer to as \emph{successive convexification}. Assuming that the solution process is initialized properly and that the original problem admits a feasible solution, then the resulting converged solution will satisfy the dynamics, kinematics, and constraints imposed in Problem~\ref{problem:non_convex}.

The remainder of this section is organized as follows. Section~\ref{sec:linearization} outlines the process of linearizing Problem~\ref{problem:non_convex}. Section~\ref{sec:discretization} outlines the process of discretizing the linearized problem. Lastly, Section~\ref{sec:successive} summarizes the convex fixed-final-time SOCP sub-problem, and provides a detailed outline of the successive convexification algorithm.

\subsection{Linearization} \label{sec:linearization}
\subsubsection{Dynamics and Kinematics} \label{sec:linearization_dynkin}

For convenience, we define the state vector $\xx(t)\in\realD{14}$ and control vector $\uu(t)\in\realD{3}$ as follows:
\begin{subequations} \label{cvx:state_ctrl_vecs}
	\begin{align}
		\xx(t) &\definedas \mat{\m(t) & \rI^T(t) & \vI^T(t) & \qIB^T(t) & \omegaB^T(t)}^T \label{cvx:state_vec} \\
		\uu(t) &\definedas \TB(t) \label{cvx:ctrl_vec}
	\end{align}
\end{subequations}
\noindent We express the dynamics and kinematics from Eqs.~\ref{prob_des:mass_depletion}-\ref{prob_des:att_kin_dyn} as a nonlinear vector-valued function $f:\realD{14}\times\realD{3}\rightarrow\realD{14}$ of the state and control vectors defined in Eq.~\ref{cvx:state_ctrl_vecs}:
\begin{equation} \label{cvx:dynamics}
	\dd{}{t}\xx(t) = \f{\xx(t)}{\uu(t)} \definedas \mat{\mdot(t)&\rIdot^T(t)&\vIdot^T(t)&\qIBdot^T(t)&\omegaBdot^T(t)}^T
\end{equation}
\noindent Next, to address the free-final-time aspect of the problem, we would like to express Eq.~\ref{cvx:dynamics} in terms of normalized trajectory time, $\tau\in[0,1]$. To do so, we apply the chain rule to the left side of Eq.~\ref{cvx:dynamics}:
\begin{equation*}
	\dd{}{t}\xx(t) = \dd{\tau}{t}\dd{}{\tau}\xx(t)
\end{equation*}
We denote the time dilation between $\tau$ and $t$ by the dilation coefficient $\sigma$:
\begin{equation} \label{cvx:sigma_eq}
	\sigma \definedas \left(\dd{\tau}{t}\right)^{-1}
\end{equation}
Eq.~\ref{cvx:dynamics} can now be rewritten in terms of $\tau$ as follows:
\begin{equation} \label{cvx:normalized_dynamics}
	\xx'(\tau) \definedas \dd{}{\tau}\xx(\tau) = \sigma\f{\xx(\tau)}{\uu(\tau)}
\end{equation}
In order to fit the nonlinear equations embedded in Eq.~\ref{cvx:normalized_dynamics} into a convex optimization framework, we approximate the right hand side with a first-order Taylor series approximation. We define a reference trajectory comprised of state, $\xxo(\tau)$, control, $\uuo(\tau)$, and dilation coefficient, $\sso$. These quantities will be explicitly defined in Section~\ref{sec:trust_regions}. The linearized system is expressed with respect to normalized time as shown below:
\begin{subequations} \label{cvx:linearization}
\begin{align}
	\xx'(\tau) &= A(\tau)\xx(\tau)+B(\tau)\uu(\tau)+\Sigma(\tau)\sigma+\zz(\tau) \\
	& \nonumber \\
	A(\tau) &\definedas \sso\thinspace\cdot\dfdx \\
	B(\tau) &\definedas \sso\thinspace\cdot\dfdu \\
	\Sigma(\tau) &\definedas \f{\xxo(\tau)}{\uuo(\tau)} \\
	\zz(\tau) &\definedas -A(\tau)\xxo(\tau)-B(\tau)\uuo(\tau)
\end{align}
\end{subequations}

\subsubsection{State Constraints}

Per the discussion in Section~\ref{sec:state_control_constraints}, the thrust magnitude lower bound constraint is the only remaining source of non-convexity in Problem~\ref{problem:non_convex}. As such, we proceed to linearize this constraint. First, we define the function $g:\realD{3}\rightarrow\realD{}$ as follows:
\begin{equation} \label{cvx:g_nonconvex}
	\g{\uu(\tau)} \definedas \Tmin-\twonorm{\uu(\tau)} \leq 0
\end{equation}
Then, replacing the left hand side of Eq.~\ref{cvx:g_nonconvex} with a first-order Taylor series approximation and simplifying, we obtain the following:
\begin{subequations} \label{cvx:g_linearization}
\begin{align}
	\Tmin &\leq B_g(\tau)\uu(\tau) \\
	B_g(\tau) &\definedas \frac{\uuo^T(\tau)}{\twonorm{\uuo(\tau)}}
\end{align}
\end{subequations}

\subsection{Discretization} \label{sec:discretization}

To cast the continuous-time optimal control problem into a finite-dimensional parameter optimization problem, we discretize the trajectory into $K$ evenly distributed discretization points. We perform the discretization with respect to normalized trajectory time, $\tau$. For convenience, we define the following two sets:
\begin{align*}
	\setK &\definedas \{0 , 1 , \dots, K-2, K-1\} \\
	\setKm &\definedas \{0 , 1, \dots, K-3, K-2 \}
\end{align*}
Since normalized trajectory time is defined on the interval $\tau\in[0,1]$, we define the time at index $k$ as:
\begin{equation*}
	\tauk \definedas \left(\frac{k}{K-1}\right) \thinspace,\quad\forall\thinspace k\in\setK
\end{equation*}
\noindent To preserve more feasibility, we assume a first-order-hold on the control over each time step. Thus, over the interval $\tau\in[\tauk,\taukp]$, we can express $\uu(\tau)$  in terms of $\uuk \definedas \uu(\tauk)$ and $\uukp \definedas \uu(\taukp)$ as follows:
\vspace{0.2cm}
\begin{subequations} \label{cvx:foh}
\begin{align}
	\uu(\tau) \definedas \lambdal(\tau)\uuk&+\lambdar(\tau)\uukp \thinspace,\quad\tau\in[\tauk,\taukp],\thinspace\forall\thinspace k\in\setKm \\
	& \nonumber \\
	\lambdal(\tau) &\definedas \left(\frac{\taukp-\tau}{\taukp-\tauk}\right) \\
	\lambdar(\tau) &\definedas \left(\frac{\tau-\tauk}{\taukp-\tauk}\right)
\end{align}
\end{subequations}
\noindent We use $\stm{\taukp}{\tauk}$ to denote the state transition matrix that describes the zero-input evolution from $\xxk\definedas\xx(\tauk)$ to $\xxkp\definedas\xx(\taukp)$. The state transition matrix is governed by the following differential equation:
\begin{equation} \label{cvx:stm}
	\frac{d}{d\tau}\stm{\tau}{\tauk} = A(\tau)\stm{\tau}{\tauk}\thinspace, \quad \stm{\tauk}{\tauk} = I, \quad\forall k\in\setKm
\end{equation}
\noindent Using Eqs.~\ref{cvx:foh} and \ref{cvx:stm}, we express the discrete-time dynamics that relate $\xxk$ to $\xxkp$ as follows:
\vspace{0.1cm}
\begin{subequations} \label{cvx:discretization}
\begin{align}
	\xxkp &= \AD\xxk+\BD\uuk+\CD\uukp+\SD\sigma+\zD \thinspace,\quad\forall\thinspace k\in\setKm \label{cvx:discrete_dyn} \\ 
	& \nonumber \\
	\AD & \definedas \stm{\taukp}{\tauk} \label{cvx:discretization_begin} \\
	\BD &\definedas \int_{\tauk}^{\taukp}{\stm{\taukp}{\xi} B(\xi)\lambdal(\xi) d\xi} \\
	\CD &\definedas \int_{\tauk}^{\taukp}{\stm{\taukp}{\xi} B(\xi)\lambdar(\xi) d\xi} \\
	\SD &\definedas \int_{\tauk}^{\taukp}{\stm{\taukp}{\xi} \Sigma(\xi) d\xi} \\
	\zD &\definedas \int_{\tauk}^{\taukp}{\stm{\taukp}{\xi} \zz(\xi) d\xi} \label{cvx:discretization_end} 
\end{align}
\end{subequations}
\noindent Lastly, the remaining state and control constraints in Problem~\ref{problem:non_convex} are enforced for each $\tauk$, $\forall\thinspace k\in\setK$.

\clearpage
\subsection{Successive Convexification} \label{sec:successive}

In this section we outline the successive convexification process, in which we indirectly solve a non-convex optimization problem by instead solving a sequence of related convex sub-problems. We denote the $i^{th}$ iterate in this sequence by the superscript $i$. Before we present the finalized SOCP sub-problem and successive convexification algorithm, we introduce two more important modifications: (1) trust regions, and (2) virtual control.

\subsubsection{Trust Regions} \label{sec:trust_regions}

In order for successive convexification to work, we must ensure that the problem remains bounded and feasible throughout the convergence process. The former issue arises when the linearization of an iterate results in constraints that admit an unbounded cost, whereas the latter issue will be discussed in Section~\ref{sec:dyn_rel}. We mitigate the possibility of unbounded solutions in each sub-problem by augmenting the cost function with terms that serve as soft trust regions defined around the previous iterate. To do so, we first define the relative quantities shown below:
\begin{subequations} \label{cvx:TR1}
	\begin{align}
		\delta\bvec{x}_k^i \definedas \xxk^i&-\bvec{x}_k^{i-1} \thinspace,\quad\forall\thinspace k\in\setK \\
		\delta\bvec{u}_k^i \definedas \uuk^i&-\bvec{u}_k^{i-1} \thinspace,\quad\forall\thinspace k\in\setK \\
		\delta\sigma^i \definedas \sigma^i&-\sigma^{i-1}
	\end{align}
\end{subequations}
Then, defining $\Deltabar^i\in\realp^K$ and $\Delta_{\sigma}^i\in\realp$, we impose the following constraints:
\begin{subequations} \label{cvx:TR2}
	\begin{align}
		\delta\bvec{x}_k^i \cdot \delta{\bvec{x}_k^i} + \delta\bvec{u}^i_k \cdot \delta\bvec{u}_k^i &\leq \bvec{e}_k\cdot\Deltabar^i \thinspace,\quad\forall\thinspace k\in\setK \\
		\delta\sigma^i \cdot \delta\sigma^i &\leq \Delta^i_\sigma
\end{align}
\end{subequations}
Lastly, defining the weight terms $\wdelta\in\realpp$ and $\wdeltasigma\in\realpp$, we augment the cost function from Problem~\ref{problem:non_convex} with the $c_{\Delta}^i\in\realpp$, defined as:
\begin{equation}\label{cvx:ctr}
	c_{\Delta}^i \definedas \wdelta\twonorm{\Deltabar^i}+\wdeltasigma\onenorm{\Delta_{\sigma}^i}
\end{equation}
Note that the trust regions are centered at $\xxk^{i-1}$, $\uuk^{i-1}$, and $\sigma^{i-1}$. Thus, when we perform the numerical integration of Eqs.~\ref{cvx:stm}, \ref{cvx:discretization_begin}-\ref{cvx:discretization_end} over the interval $\tau\in(\tauk,\taukp]$, we elect to evaluate the linearization in Eq.~\ref{cvx:linearization} about the nonlinear trajectory beginning at $\xxk^{i-1}$, and generated using the control $\uu(\tau)$ defined in Eq.~\ref{cvx:foh}. Doing so for each $k\in\setKm$ explicitly defines the linearization path $\xxo(\tau)$, $\uuo(\tau)$, and $\sso$ introduced in Section~\ref{sec:linearization_dynkin}.

\subsubsection{Virtual Control} \label{sec:dyn_rel}

The second issue alluded to in the previous section is referred to as \emph{artificial infeasibility}. Artificial infeasibility can be encountered during the convergence process when the linearization is not favorable for feasibility. For example, if the problem is linearized about an unrealistically short time-of-flight, the linearized equations will likely not admit a feasible solution. Artificial infeasibility is encountered very frequently during the first few iterations of a typical successive convexification sequence, and is largely due to initiating the process with a poor initial guess. To mitigate this commonly encountered condition, we introduce a virtual control term, $\nuk^i\in\realD{14}$, and add it to Eq.~\ref{cvx:discrete_dyn}. Consequently, the dynamics we encode in the problem are given by:
\begin{equation} \label{cvx:discrete_dyn_relaxed}
	\xxkp^i = \AD^i\xxk^i+\BD^i\uuk+\CD^i\uukp^i+\SD^i\sigma^i+\zD^i+\nuk^i \thinspace,\quad\forall\thinspace k\in\setKm
\end{equation}
For notational convenience, we concatenate the $\nuk^i$ vectors into a larger vector $\nubar^i\in\realD{14(K-1)}$:
\begin{align*}
	\nubar^i &\definedas \mat{{\bvec{\nu}_0^i}^T&\cdots&{\bvec{\nu}^i_{K-2}}^T}^T
\end{align*}
Lastly, defining the weight term $\wnu\in\realpp$, we augment the already modified cost function from Problem~\ref{problem:non_convex} with $c_{\nu}^i\in\realp$, defined as:
\begin{equation} \label{cvx:cdr}
	c_{\nu}^i \definedas \wnu\onenorm{\nubar^i}
\end{equation}
The intuition behind this modification is that by selecting a large value for $\wnu$, $\nuk^i$ serves as a heavily penalized auxiliary control variable that acts when necessary to prevent infeasibility. A converged solution with a negligible $c_{\nu}^i$ component of the cost indicates that the converged solution is dynamically feasible.

\subsubsection{Convex Sub-Problem}

We are now ready to summarize the convex sub-problem that will be solved repeatedly by the successive convexification algorithm. By virtue of the time normalization introduced in Section~\ref{sec:linearization_dynkin}, this problem can be viewed as a fixed-final-time optimization problem due to the fact that the final normalized time is always equal to unity. As a consequence, $\tf$ in Problem~\ref{problem:non_convex} is substituted with $\sigma_i$. Lastly, for convenience, we rewrite Eq.~\ref{prob_des:tilt} in more compact form:
\begin{subequations} \label{cvx:tilt}
	\begin{align}
		\cos\maxtilt &\leq 1-2\twonorm{H_q \qIBk^i}^2 \\
		H_q &\definedas \mat{0&0&1&0 \\ 0&0&0&1}
	\end{align}
\end{subequations}
The summary of the convex sub-problem is provided in Problem~\ref{problem:convex}.

\boxing{ht!}{problem}{problem:convex}{16cm}{Convex Discrete-Time Fixed-Final-Time Problem}{
	\begin{tabular}{llll}
		&&& \\
		\underline{Cost Function}: &\null\quad\null&\null\quad\null& \\
		&& $\underset{\sigma^i,\uuk^i}{\text{minimize}}\quad \sigma^i+\wnu\onenorm{\nubar^i}+\wdelta\twonorm{\Deltabar^i}+\wdeltasigma\onenorm{\Delta_{\sigma}}$ & (See Eqs.~\ref{cvx:ctr},\ref{cvx:cdr}) \\
		&&& \\
		&& subject to: &\\
		\underline{Boundary Conditions}: &&& \\
		&& $ {\begin{aligned}
					\mki^i &= \mi && \\ 
					\rIki^i &= \rIi & \rIkf^i &= \rIf \\
					\vIki^i &= \vIi & \vIkf^i &= \vIf \\
					&& \qIBkf^i &= \qIBf \\
					\omegaBki^i &= \omegaBi & \omegaBkf^i &= \omegaBf \\
					&& \ee\cdot\uukf^i &= \en\cdot\uukf^i = 0
				\end{aligned}} $ & \\
		\underline{Dynamics}: &&& \\
		&& $ \xxkp^i = \AD^i\xxk^i+\BD^i\uuk^i+\CD^i\uukp^i+\SD^i\sigma^i+\zD^i+\nuk^i $ & ${\begin{aligned}
		(&\text{See Eqs.~\ref{cvx:linearization},}\\
		&\text{\ref{cvx:foh},\ref{cvx:stm},\ref{cvx:discretization},\ref{cvx:discrete_dyn_relaxed})}\end{aligned}}$ \\
		\underline{State Constraints}: &&& \\
		&& $ {\begin{aligned}
					\mf &\leq \mk^i \\
					\tan{\glideslope}\twonorm{H_{23}\rIk^i} &\leq \eu\cdot\rIk^i \\
					\cos{\maxtilt} &\leq 1-2\twonorm{H_q\qIBk^i}^2 \\
					\twonorm{\omegaBk^i} &\leq \omegamax \\
				\end{aligned}} $ & $ {\begin{aligned}
				&\text{(See Eq.~\ref{prob_des:mdry})} \\
				&\text{(See Eq.~\ref{prob_des:gs_cone})} \\
				&\text{(See Eq.~\ref{cvx:tilt})} \\
				&\text{(See Eq.~\ref{prob_des:ang_rate})}
				\end{aligned}} $ \\
		\underline{Control Constraints}: &&& \\
		&& $ {\begin{aligned}
					\Tmin &\leq B_g(\tauk)\uuk^i \\
					\twonorm{\uuk^i} &\leq \Tmax \\
					\cos{\maxgimbalangle} \twonorm{\uuk^i} &\leq \eu\cdot\uuk
				\end{aligned}} $ & $ {\begin{aligned}
				&\text{(See Eq.~\ref{cvx:g_linearization})} \\
				&\text{(See Eq.~\ref{prob_des:thrust_mag})} \\
				&\text{(See Eq.~\ref{prob_des:gimbal_ang})}
				\end{aligned}} $ \\
		\underline{Trust Regions}: &&& \\
		&& $ {\begin{aligned}
					\delta\bvec{x}_k^i \cdot \delta\bvec{x}_k^i + \delta\bvec{u}_k^i \cdot \delta\bvec{u}_k^i &\leq \Delta_k^i \\
					\onenorm{\delta\sigma^i} &\leq \Delta_{\sigma}^i
				\end{aligned}} $ & (See Eqs.~\ref{cvx:TR1}-\ref{cvx:TR2}) \\
		&&& \\
	\end{tabular}
}

\subsubsection{Algorithm}

In this section we summarize the successive convexification algorithm. We define two tolerances, $\Delta_{tol}\in\realpp$ and $\nu_{tol}\in\realpp$ that will be used to define the exit condition. The algorithm is shown in Algorithm~\ref{algorithm:scvx}.


\boxing{ht!}{algorithm}{algorithm:scvx}{16cm}{Successive Convexification}{
	\begin{tabbing}
		\hspace{1em} \= \hspace{1em} \= \hspace{1cm} \= \hspace{0.5cm} \= \hspace{1cm} \= \hspace{1cm} \= \kill
		\kill \\
		\>{\underline{Initialization}:} \\
			\\
			\>\> for $\forall\thinspace k\in\setK$:
			\\
			\>\>\> define: $\displaystyle{\quad\alpha_1\definedas \frac{K-k}{K} \quad\quad\alpha_2\definedas \frac{k}{K}}$  \\
			\\
			\>\>\>$ \begin{aligned}
					\mk^0 &= \alpha_1\mi + \alpha_2\mf     &  \\
					\rIk^0 &= \alpha_1\rIi & \\
					\vIk^0 &= \alpha_1\vIi + \alpha_2\vIf  & \\
					\qIBk^0 &= \mat{1&0&0&0}^T & \\
					\omegaBk^0 &= \bvec{0} &
				\end{aligned} $ \hspace{8.2cm} \\
			\\
			\>\>\> set: $\xxk^0 = \mat{\mk & \rIk^T & \vIk^T & \qIBk^T & \omegaBk^T}^T$, $\uuk^0 = -\mk^0\gI$, and $\sigma^0 = t_{f,guess}$ \\
			\>\> end \\
			\\
			\>\> for $\forall\thinspace k\in\setK$: \\
			\>\>\> compute: $\AD^0$, $\BD^0$, $\CD^0$, $\SD^0$, and $\zD^0$ using Eqs.~\ref{cvx:linearization}, \ref{cvx:foh}-\ref{cvx:discretization}, \ref{cvx:discrete_dyn_relaxed} \\
			\>\> end \\
			\\
		\>{\underline{Successive Convexification Loop}:} \\
			\\
			\>\> for $i\in\{1,N_{iter,max}\}$ \\
			\\
			\>\>\> (1) solve Problem~\ref{problem:convex} using $\xxk^{(i-1)}$, $\uuk^{(i-1)}$, $\sigma^{(i-1)}$, $\AD^{(i-1)}$, $\BD^{(i-1)}$, $\CD^{(i-1)}$, $\SD^{(i-1)}$, \and $\zD^{(i-1)}$ \\[1em]
			\>\>\> (2) store newly computed $\xxk^i$, $\uuk^{i}$, and $\sigma^i$ \\[1em]
			\>\>\> (3) \> if $(\twonorm{\Deltabar^{i}} \leq \Delta_{tol}) \text{ and } (\onenorm{\nubar^i} \leq \nu_{tol})$ \\[1em]
			\>\>\>\>\> {\bf exit loop} \\[1em]
			\>\>\>\> else \\[1em]
			\>\>\>\>\> (3.1) compute new $\AD^{i}$, $\BD^{i}$, $\CD^i$, $\SD^i$, and $\zD^{i}$ using Eqs.~\ref{cvx:linearization}, \ref{cvx:foh}-\ref{cvx:discretization}, \ref{cvx:discrete_dyn_relaxed} \\[1em]
			\>\>\>\>\> (3.2) increment $i$ \\[1em]
			\>\>\>\>\> (3.3) return to step 1\\[1em]
			\>\>\>\> end \\
			\\
			\>\> end \\
		\\
	\end{tabbing}
}


\clearpage
\section{Simulation Results} \label{sec:results}

In this section we present simulation results for the proposed algorithm. These results were generated in MATLAB using CVX\cite{cvx} and SDPT3\cite{sdpt3_pap}. Two scenarios were considered: a two-dimensional in-plane maneuver, and a three-dimensional out-of-plane maneuver. While the first scenario could be equivalently modeled using a simplified formulation with two translational states and one rotational state, we emphasize that the results presented here were generated using the full 6-DoF algorithm. For the purpose of illustration, the simulation results were generated using notional non-dimensional quantities, and thus are not intended to match a real-world system. Tables~\ref{table1} and \ref{table2} contain the simulation parameters, algorithm parameters, and boundary conditions common to both scenarios. Note that in both cases, the algorithm used 50 time discretization points, and was limited to a maximum of 15 iterations.

\begin{figure*}[h!]
  \vspace{0.5cm}
  \begin{multicols}{2}
    \begin{table}[H]
      \begin{center}
        \caption{Simulation Parameters}
        \begin{tabular}{lll}
          \hhline{===}
          Parameter          & Value          			& Units \\      
          \hline
          $\gI$                & $-\eu$             & $\UL/\UT^2$ \\
					$\mi$                & $2.00$             & $\UM$ \\
					$\mf$                & $1.00$             & $\UM$ \\
					$\Tmin$              & $0.30$             & $\UM\UL/\UT^2$ \\
					$\Tmax$              & $5.00$             & $\UM\UL/\UT^2$ \\
					$\maxgimbalangle$    & $20$               & $\degree$ \\
					$\maxtilt$           & $90$               & $\degree$ \\
					$\glideslope$        & $20$               & $\degree$ \\
					$\omegamax$          & $60$               & $\degree/\UT$ \\
					$\JB$                & $1e^{-2}\cdot I_{3\times3}$   & $\UM\UL^2$ \\
					$\rTB$               & $-1e^{-2}\cdot\eu$ & $\UL$ \\
          \hhline{===}
          \label{table1}
        \end{tabular}
      \end{center}
    \end{table}
    
    \begin{table}[H]
      \begin{center}
        \caption{B.C.'s and Algorithm Parameters}
        \begin{tabular}{lll}
          \hhline{===}
          Parameter                 & Value                   & Units \\
          \hline
          $\wnu$                    & $1e^{5}$                & - \\
					$\wdelta$                 & $1e^{-3}$               & - \\
					$\wdeltasigma$            & $1e^{-1}$               & - \\
					$\nu_{tol}$               & $1e^{-10}$              & - \\
					$\Delta_{tol}$            & $1e^{-3}$               & - \\
					$N_{iter,max}$            & $15$                    & - \\
					$K$                       & $50$                    & - \\
					$\rIi$                    & $\mat{4&4&0}^T$         & $\UL$ \\
					$\vIf$                    & $-1e^{-1}\cdot\eu$      & $\UL/\UT$ \\
					$\omegaBi$                & $\bvec{0}$              & $\degree/\UT$ \\
					$\qIBf$                   & $\mat{1&0&0&0}^T$       & - \\
					\hhline{===}
          \label{table2}
        \end{tabular}
      \end{center}
    \end{table}
  \end{multicols}
  \vspace{-1.25cm}
\end{figure*}

\subsection{2-D In-Plane Example}

\begin{figure}[H]
	\begin{centering}
		\includegraphics[width=14cm]{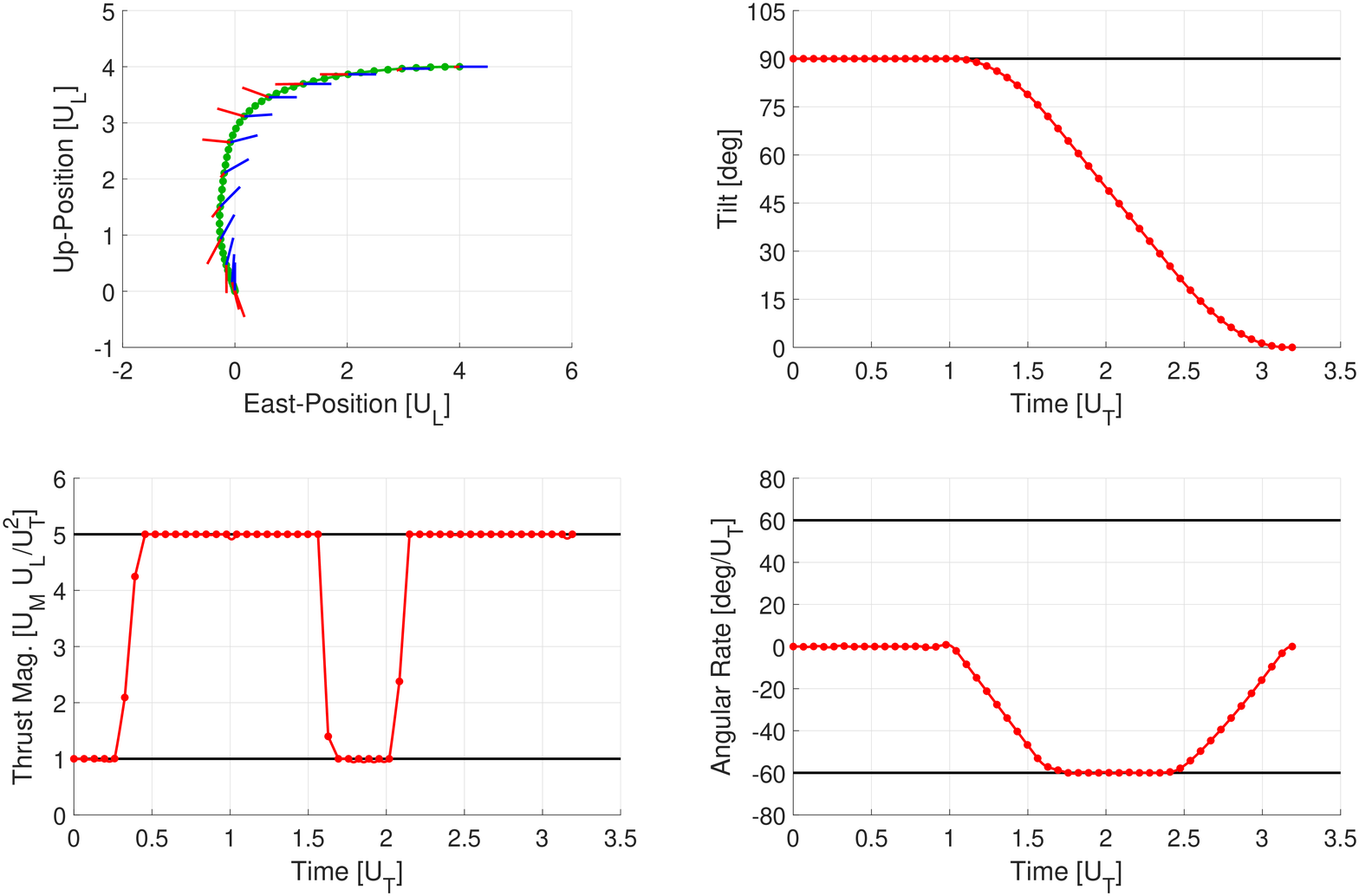}
		\vspace{0.1cm}
		\caption{2-D In-Plane Maneuver - (a) Top Left: The trajectory in the Up-East plane. The green line represents the trajectory of the vehicle, while the dots along the trajectory represent the time discretization points. The blue vectors represent the X-axis of the vehicle, whereas the red lines represent the thrust-plume direction (shown for every forth discretization point). (b) Bottom Left: The thrust magnitude profile versus time. The top and bottom black lines represent the maximum and minimum allowable thrust magnitudes, respectively. (c) Top Right: The tilt profile versus time. The black line represents the maximum allowable tilt angle. (d) Bottom Right: The angular velocity profile versus time. The black lines represent the maximum allowable angular rate.}
		\label{figure1}
	\end{centering}
\end{figure}

In this example, the initial conditions were chosen such that the vehicle remained in the Up-East plane. The trajectory was constrained to an initial horizontal velocity of 4 [$\UL/\UT$] to the west. Figure~\ref{figure1} shows the trajectory, the tilt angle, the thrust magnitude, and the angular velocity profiles. Figure~\ref{figure1}a shows the thrust and pointing vectors marked at every fourth discretization point. The trajectory is seen to ride the maximum tilt angle limit during the first half of the trajectory, and the angular velocity limit during most of the second half. Throughout, both the minimum and maximum thrust magnitude bounds are activated, as one would expect for a minimum-time trajectory.

In order to test the algorithm's robustness to the time-of-flight initial guess, the algorithm was initialized using ten different time-of-flight guesses, ranging from $1.0$ [$\UT$] to $10.0$  [$\UT$], in $1.0$ [$\UT$] increments. All ten initializations generated the same converged trajectory, yielding time-of-flights within $0.01$ [$\UT$] of each other. The red lines in Figure~\ref{figure4} show the time-of-flight as a function of iteration number for each of the ten time-of-flight initializations. In this example, convergence was obtained by the sixth iteration. These results suggest that the algorithm converged to the same solution despite the non-convexity of the problem, and despite the simple initialization routine specified in Algorithm~\ref{algorithm:scvx}.

\subsection{3-D Out-of-Plane Example}

\begin{figure}[H]
	\begin{centering}
		\includegraphics[width=14cm]{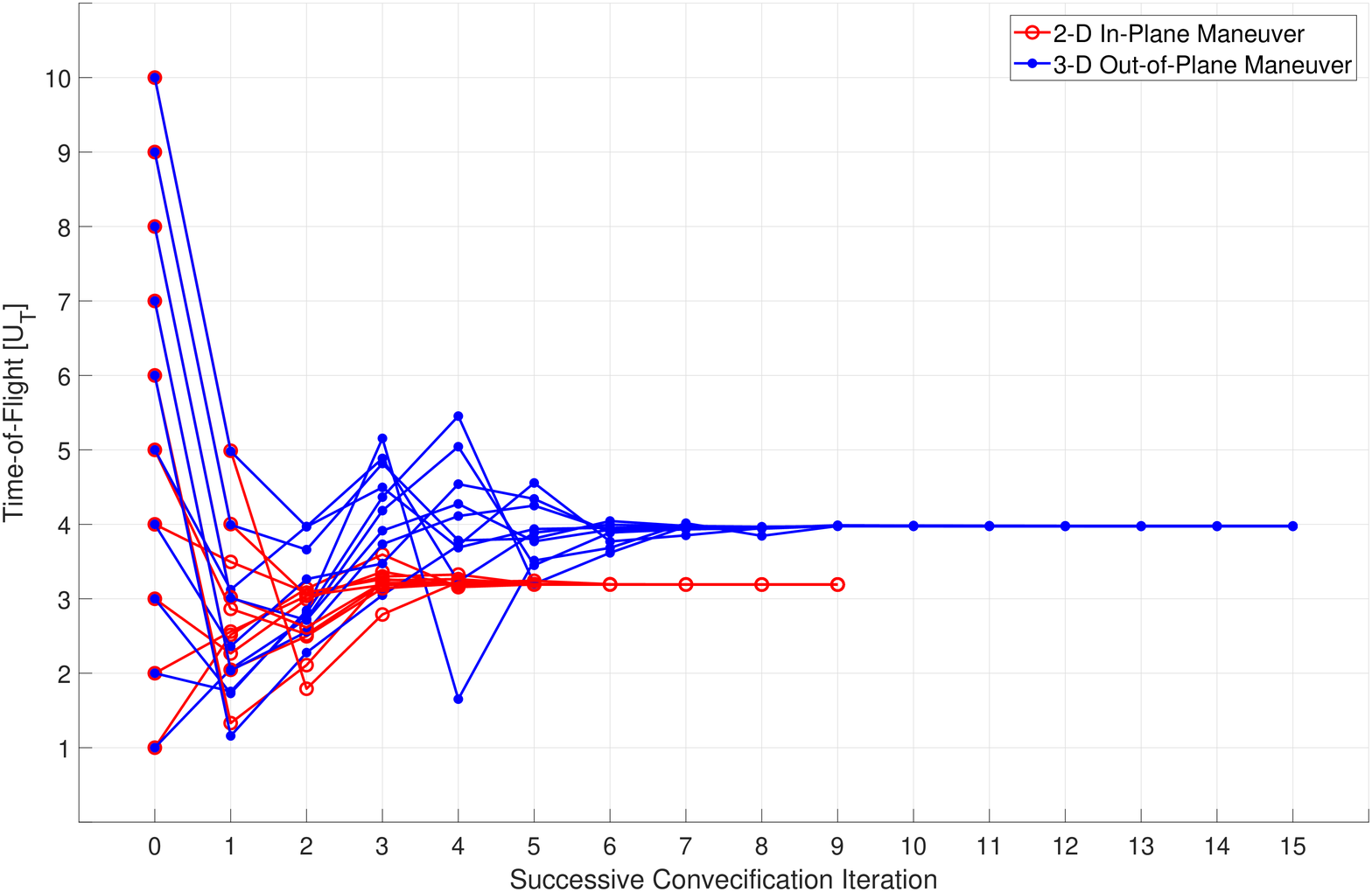}
		\vspace{0.25cm}
		\caption{Time-of-flight as a function of successive convexification iteration. Each line represents the convergence history resulting from a different time-of-flight initialization. The red and blue lines represent the 2-D and 3-D convergence histories, respectively. The fact that all lines of a given color converge to the same time-of-flight indicates that the algorithm found the same minimum objective value.}
		\label{figure4}
	\end{centering}
\end{figure}

\begin{figure}[H]
	\begin{centering}
		\includegraphics[width=14cm]{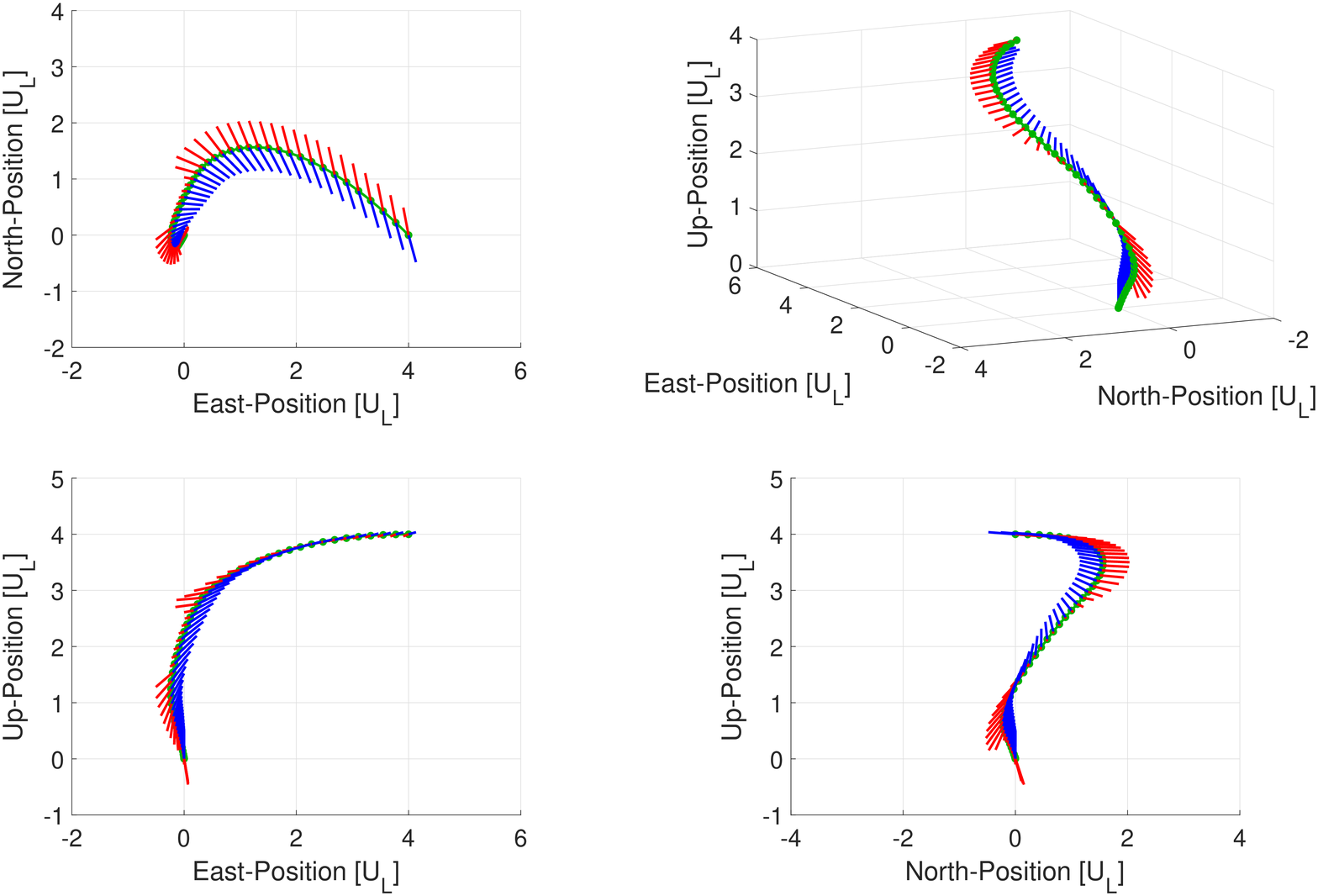}
		\caption{3-D Out-of-Plane Maneuver - (a) Top Left: North-East view of the trajectory. (b) Bottom Left: Up-East view of the trajectory. (c) Top Right: Perspective view of the trajectory. (d) Bottom Right: Up-North view of the trajectory.}
		\label{figure2}
	\end{centering}
\end{figure}

\begin{figure}[H]
	\begin{centering}
		\includegraphics[width=14cm]{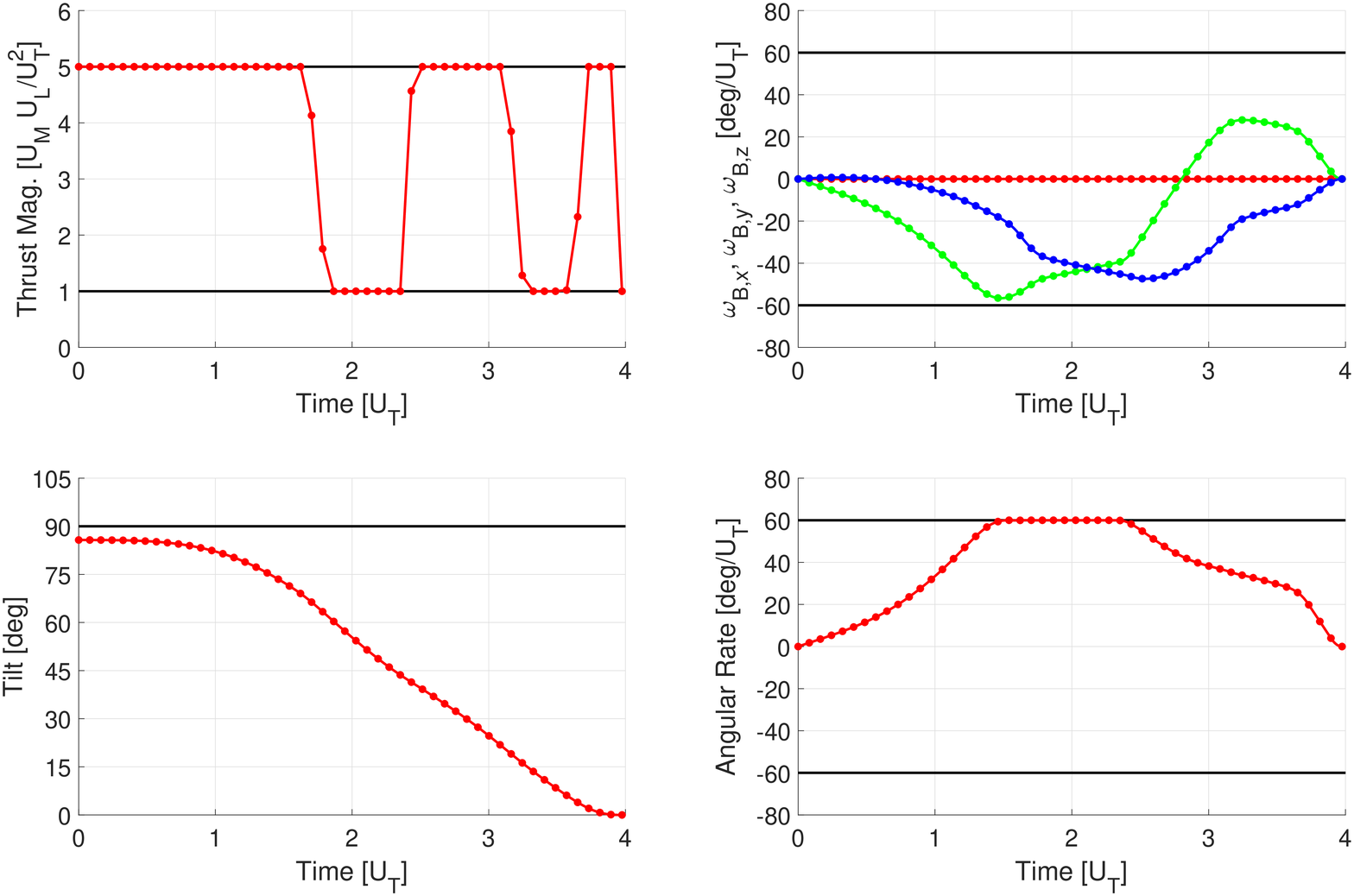}
		\caption{3-D Out-of-Plane Maneuver - (a) Top Left: Thrust magnitude profile versus time. (b) Bottom Left: Tilt angle profile versus time. (c) Top Right: Angular rate components versus time. The x, y, and z components are represented by the red, green, and blue lines, respectively. (d) Bottom Right: Angular rate magnitude versus time.}
		\label{figure3}
	\end{centering}
\end{figure}

Figure~\ref{figure2} and \ref{figure3} show results for the 3-D out-of-plane example. Despite the more complex nature of the trajectory, the algorithm performed nearly identically when compared to the 2-D case. Here, bang-coast-bang features and angular rate saturation similar to those observed in the 2-D case are apparent in the thrust magnitude and angular rate profiles, respectively. As seen from the blue lines in Figure~\ref{figure4}, the 3-D case also converged to the same time-of-flight in all ten initializations, achieving convergence after nine iterations.

\section{Conclusion and Future Work} \label{sec:conclusion}

In this paper, we presented a successive convexification formulation to solve the minimum-time 6-DoF Mars rocket powered landing problem. To solve this nonlinear, non-convex optimal control problem, we employed the following procedure. First, we initialized the process with an easily generated, dynamically inconsistent trajectory. Second, for each iterate, (a) the problem was linearized and discretized by computing the state and control matrices for the corresponding linear-time-varying system, (b) non-convex constraints were linearized, and (c) the resulting convex SOCP sub-problem was solved to full optimality. Third, the process was repeated until the changes between iterates was sufficiently small, and the virtual control term subsided.

Simulations were conducted in order to exercise the proposed algorithm. The results presented showed that the algorithm can indeed be initialized with a simple trajectory and a poor time-of-flight guess. Despite the complexity of the problem, this method was shown to converge in less than 15 successive convexification iterations.

In future work, we plan to explore the convergence properties of the successive convexification technique outlined in this paper, present real-time timing statistics, and incorporate simple aerodynamic effects into the 6-DoF problem.

\section*{Acknowledgments}

Support for studying the convergence properties of the successive convexification framework was provided by the Office of Naval Research grants N00014-16-1-2877 and N00014-16-1-3144.

\bibliography{bibliography}

\end{document}